\documentclass[12pt]{amsart}
\usepackage{amssymb}
\usepackage{graphics}

\def\@abssec#1{\vspace{.05in}\footnotesize \parindent .2in
{\bf #1. }\ignorespaces}

\newtheorem*{thm}{Theorem}
\newtheorem*{lemma}{Lemma}

\def\th{\theta}
\def\om{\omega}
\def\R{\Bbb R}
\def\cP{\mathcal P}

\renewcommand{\le}{\leqslant}
\renewcommand{\ge}{\geqslant}

\begin{document}

\title[Well-posedness for quasi-geostrophic
equation]{Global well-posedness for the critical $2D$ dissipative
quasi-geostrophic equation}
\author{A. Kiselev}
\thanks{Department of
Mathematics, University of Wisconsin, Madison, WI 53706; e-mail:
kiselev@math.wisc.edu}
\author{F. Nazarov}
\thanks{Department of
Mathematics, Michigan State University, East Lansing, MI 48824;
e-mail: fedja@math.msu.edu}
\author{A. Volberg}
\thanks{Department of
Mathematics, Michigan State University, East Lansing, MI 48824;
e-mail: volberg@math.msu.edu}


\begin{abstract}
We give an elementary proof of the global well-posedness for the
critical $2D$ dissipative quasi-geostrophic equation. The argument
is based on a non-local maximum principle involving appropriate
moduli of continuity.
\end{abstract}
\maketitle

\section{Introduction and main results}

The $2D$ quasi-geostrophic equation attracted quite a lot of
attention lately from various authors. Mainly it is due to the
fact that it is the simplest evolutionary fluid dynamics equation
for which the problem of existence of smooth global solutions
remains unsolved. In this paper we will consider the so-called
dissipative quasi-geostrophic equation
\begin{equation}
\left\{
\aligned
&\th_t=u\cdot\nabla\th -(-\Delta)^\alpha\th
\\ \nonumber
& u=(u_1,u_2)=(-R_2\th, R_1\th)
\endaligned
\right.
\end{equation}
where $\th\,:\,\Bbb R^2\to \Bbb R$ is a scalar function, $R_1$ and
$R_2$ are the usual Riesz transforms in $\Bbb R^2$ and $\alpha>0$.
It is well known (see \cite{CW,R}) that for $\alpha>\frac12$ (the
so-called subcritical case), the initial value problem
$\th(x,0)=\th_0(x)$ with $C^\infty$-smooth periodic initial data
$\th_0$ has a global $C^\infty$ solution.

For $\alpha=\frac12,$ this equation arises in geophysical studies
of strongly rotating fluid flows (see e.g. \cite{C} for further
references).
Therefore, a significant amount of research focused specifically
on the critical $\alpha=\frac12$ case. In particular, Constantin,
Cordoba, and Wu in \cite{CCW} showed that the global smooth
solution exists provided that $\|\th_0\|_\infty$ is small enough.
Cordoba and Cordoba \cite{CC} proved that the viscosity solutions
are smooth on time intervals $t \leq T_1$ and $t \geq T_2.$ The
aim of this paper is
to demonstrate that, in the critical case, smooth global solutions
exist for any $C^\infty$ periodic initial data $\th_0,$ with no
additional qualifications or assumptions. What happens in the
supercritical case $0\le\alpha<\frac12$ remains an open question.

The main idea of our proof is quite simple: we will construct a
special family of moduli of continuity that are preserved by the
dissipative evolution, which will allow us to get an a priori
estimate for $\|\nabla\th\|_\infty$ independent of time. More
precisely, we will prove the following theorem.
\begin{thm}\label{thm1}
The quasi-geostrophic equation with periodic smooth initial data
$\theta_0(x)$ has a unique global smooth solution. Moreover, the
following estimate holds:
\begin{equation}\label{mainest}
\|\nabla\th\|_\infty\le C\|\nabla\th_0\|_\infty\exp\exp
\{C\|\th_0\|_\infty\}\,.
\end{equation}
\end{thm}
At this moment we do not know how sharp the upper bound
\eqref{mainest} is. On the other hand, any a-priori bound for
$\|\nabla\th\|_\infty$ is sufficient for the proof of
well-posedness. Indeed, local existence and regularity results
then allow to extend the unique smooth solution indefinitely.
For the critical and supercritical quasi-geostrophic equation,
such results can be found for example in \cite{Wu} (Theorems 3.1
and 3.3). Hence, the rest of the paper is devoted to the proof of
\eqref{mainest}.

This paper is built upon the ideas discovered in a related work on
the dissipative Burgers equation \cite{KNS}.

\section{Moduli of continuity}

Let us remind the reader that a modulus of continuity is just an arbitrary
increasing continuous concave function $\om\,:\,[0,+\infty)\to[0,+\infty)$
such that $\om(0)=0$. Also, we say that a function $f\,:\,\Bbb R^n\to\Bbb
R^m$ has modulus of continuity $\om$ if $|f(x)-f(y)|\le\om(|x-y|)$ for all
$x,y\in\Bbb R^n$.

Singular integral operators like Riesz transforms do not preserve
moduli of continuity in general but they do not spoil them too
much either. More precisely, we have
\begin{lemma}\label{lemma1} If the function $\th$ has
modulus of continuity $\om$, then $u=(-R_2\th, R_1\th)$ has
modulus of continuity
\[ \Omega(\xi)=A\left(\int_0^\xi \frac{\om(\eta)}{\eta}\,d\eta+
\xi\int_\xi^\infty \frac{\om(\eta)}{\eta^2}\,d\eta\right) \]
with some universal constant $A>0$.
\end{lemma}
The proof of this result is elementary but since we could not
readily locate it in the literature, we provide a sketch in the
appendix.

The flow term $u\cdot\nabla\th$ in the dissipative quasi-geostrophic equation
tends to make the modulus of continuity of $\th$ worse while the dissipation
term $(-\Delta)^\alpha\th$ tends to make it better. Our aim is to construct
some special moduli of continuity for which the dissipation term always
prevails and such that every periodic $C^\infty$-function $\th_0$ has one of
these special moduli of continuity.

Note that  the critical ($\alpha=\frac12$) equation has a simple scaling
invariance: if $\th(x,t)$ is a solution, then so is $\th(Cx,Ct)$. This means
that if we manage to find one modulus of continuity $\om$ that is preserved
by the dissipative evolution, then the whole family $\om_C(\xi)=\om(C\xi)$ of
moduli of continuity will also be preserved.

Observe now that if $\om$ is unbounded, then any $C^\infty$ periodic function
has modulus of continuity $\om_C$ if $C>0$ is sufficiently large. Also, if
the modulus of continuity $\om$ has finite derivative at $0$, it can be used
to estimate $\|\nabla\th\|_\infty$. Thus, our task reduces to constructing an
unbounded modulus of continuity with finite derivative at $0$ that is
preserved by the dissipative evolution.

From now on, we will also assume that, in addition to
unboundedness and the condition $\om'(0)<+\infty$, we have
$\lim_{\xi\to 0+}\om''(\xi)=-\infty$. Then, if a $C^\infty$
periodic function $f$ has modulus of continuity $\om$, we have
$$
\|\nabla f\|_\infty<\om'(0)\,.
$$
Indeed, take a point $x\in\R^2$ at which $\max|\nabla f|$ is attained and
consider the point $y=x+\xi e$ where $e=\frac{\nabla f}{|\nabla f|}$. Then we
must have $f(y)-f(x)\le \om(\xi) $ for all $\xi\ge 0$. But the left hand side
is at least $|\nabla f(x)|\xi-C\xi^2$ where $C=\frac 12\|\nabla^2 f\|_\infty$
while the right hand side can be represented as $\om'(0)\xi-\rho(\xi)\xi^2$
with $\rho(\xi)\to+\infty$ as $\xi\to 0+$. Thus $|\nabla f(x)|\le \om'(0)-
(\rho(\xi)-C)\xi$ for all $\xi>0$ and it remains to choose some $\xi>0$
satisfying $\rho(\xi)>C$.

\section{The breakthrough scenario}

Now assume that $\th$ has modulus of continuity $\om$ for all times $t<T$.
Then $\th$ remains $C^\infty$ smooth up to $T$ and, according to the local
regularity theorem, for a short time beyond $T$. By continuity, we see that
$\th$ must also have modulus of continuity $\om$ at the moment $T$. Suppose
that $|\th(x,T)-\th(y,T)|<\om(|x-y|)$ for all $x\ne y$. We claim that then
$\th $ has modulus of continuity $\om$ for all $t>T$ sufficiently close to
$T$. Indeed, by the remark above, at the moment $T$ we have
$\|\nabla\th\|_\infty<\om'(0)$. By continuity of derivatives, this also holds for
$t>T$ close to $T$, which immediately takes care of the inequality
$|\th(x,t)-\th(y,t)|<\om(|x-y|)$ for small $|x-y|$. Also, since $\om$ is
unbounded and $\|\th\|_\infty$ doesn't grow with time, we automatically
have $|\th(x,t)-\th(y,t)|<\om(|x-y|)$ for large $|x-y|$. The last observation
is that, due to periodicity of $\th$, it suffices to check the inequality
$|\th(x,t)-\th(y,t)|<\om(|x-y|)$ for $x$ belonging to some compact set
$K\subset \R^2$. Thus, we are left with the task to show that, if
$|\th(x,T)-\th(y,T)|<\om(|x-y|)$ for all $x\in K$, $\delta\le |x-y|\le
\delta^{-1}$ with some fixed $\delta>0$, then the same inequality remains
true for a short time beyond $T$. But this immediately follows from the
uniform continuity of $\theta$.

This implies that the only scenario in which the modulus of continuity $\om$
may be lost by $\th$ is the one in which there exists a moment $T>0$ such
that $\om$ has modulus of continuity $\om$ for all $t\in [0,T]$ and there are
two points $x\ne y$ such that $\th(x,T)-\th(y,T)=\om(|x-y|)$. We shall rule
this scenario out by showing that, in such case, the derivative
$\bigl.\frac{\partial}{\partial t}(\th(x,t)-\th(y,t))\bigr|_{t=T}$ must be
negative, which, clearly, contradicts the assumption that the modulus of
continuity $\om$ is preserved up to the time $T$.

\section{Estimate of the derivative: the flow term}

Assume that the above scenario takes place. Let $\xi=|x-y|$. Observe that
$(u\cdot \nabla \th)(x)=\bigl.\frac{d}{dh}\th(x+hu(x))\bigr|_{h=0}$ and
similarly for $y$. But
$$
\th(x+hu(x))-\th(y+hu(y))\le\om(|x-y|+h|u(x)-u(y)|)\le \om(\xi+h\Omega(\xi))
$$
where, as before,
$$
\Omega(\xi)=A\left(\int_0^\xi \frac{\om(\eta)}{\eta}\,d\eta+
\xi\int_\xi^\infty \frac{\om(\eta)}{\eta^2}\,d\eta\right)\,.
$$
Since $\th(x)-\th(y)=\om(\xi)$, we conclude that
$$
(u\cdot\nabla\th)(x)-(u\cdot\nabla\th)(y)\le \Omega(\xi)\om'(\xi)\,.
$$

\section{Estimate of the derivative: the dissipation term}

Recall that the dissipative term can be written as
$\bigl.\frac{d}{dh}\cP_h * \th\bigr|_{h=0}$ where $\cP_h$ is the
$2$-dimensional Poisson kernel. Thus, our task is to estimate
$(\cP_h*\th)(x)-(\cP_h*\th)(y)$ under the assumption that $\th$ has modulus
of continuity $\om$. Since everything is translation and rotation invariant,
we may assume that $x=(\frac{\xi}{2},0)$ and $y=(-\frac{\xi}{2},0)$.

Write
\begin{eqnarray*} (\cP_h*\th)(x)-(\cP_h*\th)(y)=\iint_{\R^2}
[\cP_h(\tfrac{\xi}{2}-\eta,\nu)-\cP_h(-\tfrac{\xi}{2}-\eta,\nu)]\th(\eta,\nu)\,
d\eta d\nu
\\
=\int_{\R}d\nu\int_0^\infty
[\cP_h(\tfrac{\xi}{2}-\eta,\nu)-\cP_h(-\tfrac{\xi}{2}-\eta,\nu)]
[\th(\eta,\nu)-\th(-\eta,\nu)]\,d\eta
\\
\le
\int_{\R}d\nu\int_0^\infty
[\cP_h(\tfrac{\xi}{2}-\eta,\nu)-\cP_h(-\tfrac{\xi}{2}-\eta,\nu)]
\om(2\eta)\,d\eta
\\
=
\int_0^\infty[P_h(\tfrac{\xi}{2}-\eta)-P_h(-\tfrac{\xi}{2}-\eta)]
\om(2\eta)\,d\eta
\\
=
\int_0^\xi P_h(\tfrac{\xi}{2}-\eta)\om(2\eta)\,d\eta+
\int_0^\infty P_h(\tfrac{\xi}{2}+\eta)[\om(2\eta+2\xi)-\om(2\eta)]\,d\eta
\end{eqnarray*}
where $P_h$ is the $1$-dimensional Poisson kernel. Here we used symmetry and
monotonicity of the Poisson kernels together with the observation that
$\int_{\R}\cP_h(\eta,\nu)\,d\nu=P_h(\eta)$. The last formula can also be
rewritten as
$$
\int_0^{\frac{\xi}{2}}P_h(\eta)[\om(\xi+2\eta)+\om(\xi-2\eta)]\,d\eta+
\int_{\frac{\xi}{2}}^\infty
P_h(\eta)[\om(2\eta+\xi)-\om(2\eta-\xi)]\,d\eta\,.
$$
Recalling that $\int_0^\infty P_h(\eta)\,d\eta=\frac{1}{2}$, we see that the
difference $(\cP_h*\th)(x)-(\cP_h*\th)(y)-\om(\xi)$ can be estimated from
above by
\begin{eqnarray*}
\int_0^{\frac{\xi}{2}}P_h(\eta)[\om(\xi+2\eta)+\om(\xi-2\eta)-2\om(\xi)]
\,d\eta
\\
+
\int_{\frac{\xi}{2}}^\infty
P_h(\eta)[\om(2\eta+\xi)-\om(2\eta-\xi)-2\om(\xi)]\,d\eta\,.
\end{eqnarray*}
Recalling the explicit formula for $P_h$,
dividing by $h$ and passing to the limit as $h\to 0+$, we finally conclude
that the contribution of the dissipative term to our derivative is bounded
from above by
\begin{eqnarray*}
\frac{1}{\pi}\int_0^{\frac{\xi}{2}}\frac{\om(\xi+2\eta)+\om(\xi-2\eta)-2\om(\xi)}{\eta^2}
\,d\eta
\\
+
\frac{1}{\pi}\int_{\frac{\xi}{2}}^\infty
\frac{\om(2\eta+\xi)-\om(2\eta-\xi)-2\om(\xi)}{\eta^2}\,d\eta\,.
\end{eqnarray*}
Note that due to concavity of $\om$, both terms are strictly negative.

\section{The explicit formula for the modulus of continuity}

We will construct our special modulus of continuity as follows. Choose two
small positive numbers $\delta>\gamma>0$ and define the continuous function $\om$ by
$$
\om(\xi)=\xi-\xi^{\frac{3}{2}}\qquad\text{when }0\le\xi\le\delta
$$
and
$$
\om'(\xi)=\frac{\gamma}{\xi(4+\log (\xi/\delta))}\qquad\text{when }\xi>\delta\,.
$$
Note that, for small $\delta$, the left derivative of $\om$ at
$\delta$ is about $1$ while the right derivative equals
$\frac{\gamma}{4\delta}<\frac{1}{4}$. So $\om$ is concave if
$\delta$ is small enough. It is clear that $\om'(0)=1$,
$\lim_{\xi\to 0+}\om''(\xi)=-\infty$ and that $\om$ is unbounded
(it grows at infinity like double logarithm). The hard part, of
course, is to show that, for this $\om$, the negative contribution
to the time derivative coming from the dissipative term prevails
over the positive contribution coming from the flow term. More
precisely, we have to check the inequality
\begin{eqnarray*}
 A\left[ \int_0^\xi\frac{\om(\eta)}{\eta}\,d\eta+
\xi\int_\xi^\infty\frac{\om(\eta)}{\eta^2}\,d\eta \right]\om'(\xi)
+
\frac{1}{\pi}\int_0^{\frac{\xi}{2}}\frac{\om(\xi+2\eta)+\om(\xi-2\eta)-2\om(\xi)}{\eta^2}
\,d\eta
\\
+ \frac{1}{\pi} \int_{\frac{\xi}{2}}^\infty
\frac{\om(2\eta+\xi)-\om(2\eta-\xi)-2\om(\xi)}{\eta^2}\,d\eta <0
\qquad\text {for all }\xi>0\,.
\end{eqnarray*}

\section{Checking the inequality: case $0\le\xi\le\delta$}

Let $0\le\xi\le\delta$. Since $\om(\eta)\le\eta$ for all $\eta\ge 0$, we
have $\int_0^\xi \frac{\om(\eta)}{\eta}\,d\eta\le\xi$ and
$\int_\xi^\delta \frac{\om(\eta)}{\eta^2}\,d\eta\le\log\frac{\delta}{\xi}$.
Now,
$$
\int_\delta^{\infty}\frac{\om(\eta)}{\eta^2}\,d\eta=
\frac{\om(\delta)}{\delta}+\gamma
\int_\delta^\infty \frac{1}{\eta^2(4+\log(\eta/\delta))}\,d\eta
\le 1+\frac{\gamma}{4\delta}<2\,.
$$
Observing that $\om'(\xi)\le 1$, we conclude that
the positive part of the left hand side is bounded by
$A\xi(3+\log\frac{\delta}{\xi})$.

To estimate the negative part, we just use the second order Taylor formula
and monotonicity of $\om''$ on $[0,\xi]$ to get the bound
$$
\frac{1}{\pi}\int_0^{\frac{\xi}{2}}\frac{\om(\xi+2\eta)+\om(\xi-2\eta)-2\om(\xi)}{\eta^2}
\,d\eta\le
\frac{1}{\pi}\xi\om''(\xi)=-\frac{3}{4\pi}\xi\xi^{-\frac{1}{2}}\,.
$$
But, obviously,
$\xi\left(A(3+\log\frac{\delta}{\xi})-\frac{3}{4\pi}\xi^{-\frac{1}{2}}\right)<0$
on $(0,\delta]$ if $\delta$ is small enough.

\section{Checking the inequality: case $\xi\ge\delta$}

In this case, we have $\om(\eta)\le\eta$ for $0\le\eta\le\delta$ and
$\om(\eta)\le \om(\xi)$ for $\delta\le \eta\le\xi$. Hence
$$
\int_0^\xi\frac{\om(\eta)}{\eta}\,d\eta\le
\delta+\om(\xi)\log\frac{\xi}{\delta}\le
\om(\xi)\left(2+\log\frac{\xi}{\delta}\right)
$$
because $\om(\xi)\ge\om(\delta)>\frac{\delta}{2}$ if $\delta$ is small
enough.

Also
$$
\int_\xi^\infty\frac{\om(\eta)}{\eta^2}\,d\eta =\frac{\om(\xi)}{\xi}+
\gamma\int_{\xi}^\infty \frac{d\eta}{\eta^2(4+\log(\eta/\delta))}\le
\frac{\om(\xi)}{\xi}+\frac{\gamma}{\xi}\le \frac{2\om(\xi)}{\xi}
$$
if $\gamma<\frac{\delta}{2}$ and $\delta$ is small enough.

Thus, the positive term on the left hand side is bounded from above by
the expression
$A\om(\xi)\left(4+\log\frac{\xi}{\delta}\right)\om'(\xi)
=A\gamma\frac{\om(\xi)}{\xi}$.

To estimate the negative term, note that, for $\xi\ge \delta$, we have
$$
\om(2\xi)\le \om(\xi)+\frac{\gamma}{4}\le \frac{3}{2}\om(\xi)
$$
under the same assumptions on $\gamma$ and $\delta$ as above. Also, due to
concavity, we have $\om(2\eta+\xi)-\om(2\eta-\xi)\le\om(2\xi)$ for all
$\eta\ge\frac{\xi}{2}$. Therefore,
$$
\frac{1}{\pi}
\int_{\frac{\xi}{2}}^\infty
\frac{\om(2\eta+\xi)-\om(2\eta-\xi)-2\om(\xi)}{\eta^2}\,d\eta
\le
-\frac{1}{2\pi}\int_{\frac{\xi}{2}}^\infty \frac{\om(\xi)}{\eta^2}\,d\eta=
-\frac{1}{\pi}\frac{\om(\xi)}{\xi}\,.
$$
But $\frac{\om(\xi)}{\xi}(A\gamma-\frac{1}{\pi})<0$ if $\gamma$ is small
enough.

\section{Concluding remarks}

Here we just want to quote (with necessary minor modifications) a
paragraph from \cite{CC}. Note that it was written just $2$ years
ago.

\it
The case $\alpha=\frac{1}{2}$ is specially relevant because the viscous term
$(-\Delta)^\frac{1}{2}\th$ models the so-called Eckmann's pumping, which has
been observed in quasi-geostrophic flows. On the other hand, several authors
have emphasized the deep analogy existing between the dissipative
quasi-geostrophic equation with $\alpha=\frac{1}{2}$ and the $3D$
incompressible Navier-Stokes equations.

\rm This paper provides an elementary treatment of the
$\alpha=\frac{1}{2}$ case. Unfortunately, the argument does not
seem to extend to the Navier-Stokes equations due to the different
structure of nonlinearity. So, while our paper resolves the global
existence and regularity question in a physically relevant model,
it also suggests that there is a significant structural difference
between the critical $2D$ quasi-geostrophic equation and $3D$
Navier-Stokes equations.

\section{Appendix}

Here we provide a sketch of the proof of the Lemma.
\begin{proof}
The Riesz transforms are singular integral operators with kernels
$K(r,\zeta) = r^{-2} \Omega(\zeta),$ where $(r,\zeta)$ are the
polar coordinates. The function $\Omega$ is smooth and $\int_{S^1}
\Omega (\zeta) d\sigma(\zeta)=0.$ Assume that the function $f$
satisfies $|f(x)-f(y)|\leq \omega(|x-y|)$ for some modulus of
continuity $\omega.$ Take any $x,y$ with $|x-y|=\xi,$ and consider
the difference
\begin{equation}\label{rt}
P.V.\int K(x-t)f(t)\,dt - P.V.\int K(y-t)f(t)\,dt
\end{equation}
with integrals understood in the principal value sense. Note that
\[ \left| P.V. \int_{|x-t| \leq 2\xi} K(x-t) f(t)\,dt \right| =
\left| P.V. \int_{|x-t| \leq 2\xi} K(x-t) (f(t)-f(x))\,dt \right|
\leq C \int_0^{2\xi} \frac{\omega(r)}{r}\,dr. \] Since $\omega$ is
concave, we have
\[ \int_0^{2\xi} \frac{\omega(r)}{r}\,dr \leq 2 \int_0^{\xi}
\frac{\omega(r)}{r}\,dr. \] A similar estimate holds for the
second integral in \eqref{rt}. Next, let $\tilde{x} =
\frac{x+y}{2}.$ Then
\begin{eqnarray*}\left| \int_{|x-t| \geq 2\xi}
K(x-t)f(t)\,dt - \int_{|y-t| \geq 2\xi} K(y-t)f(t)\,dt\right| =\\
\left| \int_{|x-t| \geq 2\xi} K(x-t)(f(t)-f(\tilde{x}))\,dt -
\int_{|y-t| \geq 2\xi} K(y-t)(f(t)-f(\tilde{x}))\,dt\right| \\
\leq \int_{|\tilde{x}-t| \geq
3\xi}|K(x-t)-K(y-t)||f(t)-f(\tilde{x})|\,dt + \\
\int_{3\xi/2 \leq |\tilde{x}-t| \leq
3\xi}(|K(x-t)|+|K(y-t)|)|f(t)-f(\tilde{x})|\,dt.
\end{eqnarray*}
Since
\[ |K(x-t)-K(y-t)| \leq C\frac{|x-y|}{|\tilde{x}-t|^3} \] when $|\tilde{x}-t| \geq
3\xi,$ the first integral is estimated by $C\xi \int_{3\xi}^\infty
\frac{\omega(r)}{r^2}\,dr.$ The second integral is estimated by
$C\omega(3\xi),$ and hence is controlled by $3C\int_0^{\xi}
\frac{\omega(r)}{r}\,dr.$
\end{proof}

{\bf Acknowledgement} Research of AK has been partially supported
by the NSF-DMS grant 0314129. Research of FN and AV has been
partially supported by the NSF-DMS grant 0501067.

\end{document}